\documentclass[12pt]{article}
\usepackage{amssymb, amsmath, amsthm, amsgen, amstext, amsbsy, amsopn}
\usepackage{amsfonts}

\newtheorem{thm}{Theorem}[section]

\newtheorem{lemma}[thm]{Lemma}

\newtheorem{cor}[thm]{Corollary}

\newtheorem{rem}[thm]{Remark}

\newtheorem{example}[thm]{Example}

\newcommand{\R}{{\mathbb{R}}}

\newcommand{\Z}{{\mathbb{Z}}}

\newcommand{\cA}{{\mathcal{A}}}

\newcommand{\cF}{{\mathcal{F}}}

\newcommand{\cH}{{\mathcal{H}}}

\newcommand{\cV}{{\mathcal{V}}}

\newcommand{\Qed}{\hfill $\Box$ \medskip}

\begin{document}
\title{Nodal inequalities on surfaces}
\author{ Leonid Polterovich and Mikhail Sodin\\
School of Mathematical Sciences\\ Tel-Aviv University}

\date{}

\maketitle

\begin{abstract}
Given a Laplace eigenfunction on a surface, we study the
distribution of its extrema on the nodal domains. It is
classically known that the absolute value of the eigenfunction is
asymptotically bounded by the 4-th root of the eigenvalue. It
turns out that the number of nodal domains where the eigenfunction
has an extremum of such order, remains bounded as the eigenvalue
tends to infinity. We also observe that certain restrictions on
the distribution of nodal extrema and a version of the Courant
nodal domain theorem are valid for a rather wide class of
functions on surfaces. These restrictions follow from a bound in
the spirit of Kronrod and Yomdin on the average number of
connected components of level sets.
\end{abstract}

\section{Introduction and main results} \label{sec-main}

Let $M$ be a compact connected surface, which in the case when
$\partial M \neq \emptyset$ is assumed to be oriented. A {\it
nodal domain} of a function $f$ on $M$ is a connected component of
the set $\{f \neq 0\}$. We write $\cA(f)$ for the collection of
nodal domains of $f$. In this note we are interested in the
distribution of nodal extrema $\displaystyle m_A :=\max_A |f|$,
where $A \in \cA(f)$.

Let $g$ be a Riemannian metric on $M$. We write $||f||$ for the
$L^2$-norm of a function $f$ on $M$ with respect to the Riemannian
area $\sigma$ on $M$. We denote by $\Delta$ the Laplace-Beltrami
operator.

Consider the space $\cF$ of all smooth functions on $M$ which, in
the case when $\partial M \neq \emptyset$, are assumed to vanish
on $\partial M$. Denote by $\cF_{\lambda}$ the subset of all
functions $f \in \cF$ with $||f||=1$, $||\Delta f|| \leq \lambda$
and, if $M$ has no boundary, $\int_M f d\sigma =0$. Note that
$\cF_{\lambda}$ is empty for $\lambda < \lambda_1(M,g)$ where
$\lambda_1$ stands for the first (Dirichlet) eigenvalue. This
follows from the variational principle for $\lambda_1$ combined
with the estimate
$$||\nabla f||^2 \leq ||f||\cdot||\Delta f|| \leq \lambda ||f||^2$$
for every $f \in \cF_{\lambda}$. Therefore we shall assume that
$\lambda \geq \lambda_1$. Let us present some examples of
functions from $\cF_\lambda$:

\medskip

\begin{example} \label{ex-1} {\rm The class $\cF_{\lambda}$ contains
the (Dirichlet) eigenfunctions $f_{\lambda}$ of the
Laplace-Beltrami operator on a compact surface, that is functions
satisfying
$$\Delta f_{\lambda} +\lambda f_{\lambda} = 0\;,\;\; f|_{\partial M} = 0\;,$$
normalized by $||f|| = 1$. Furthermore, any normalized linear
combination of eigenfunctions  with eigenvalues $\lambda_i \leq
\lambda$ belongs to $\cF_{\lambda}$.}
\end{example}

\begin{example} \label{ex-4}
{\rm Another interesting class of functions from $\cF_{\lambda}$
is given by normalized eigenfunctions of the biharmonic operator
on a surface with boundary with the ``clamped plate" boundary
conditions:
$$\Delta^2 f = \lambda^2 f,\;\; f|_{\partial M}=0,\;\nabla
f|_{\partial M}=0\;.$$}
\end{example}

\medskip
\noindent We start with the following result:

\begin{thm}\label{thm-crit} For any $f\in\cF_{\lambda}$,
\begin{equation}\label{eq-gen-11}
\displaystyle \sum_{A\in \cA} m_A \leq k_g \lambda\,,
\end{equation}
and
\begin{equation}\label{eq-gen-21}
\displaystyle \sum_{A\in \cA} m_A^2 \leq k_g \lambda\,.
\end{equation}
\end{thm}

\noindent Here and in what follows, by $k_g$ we denote positive
constants that depend only on the metric $g$. The value of these
constants may vary from line to line.

Recall that the Courant nodal domain theorem combined with the
Weyl law \cite{Ch} yields that the number of nodal domains of a
(Dirichlet) eigenfunction $f_{\lambda}$ of the Laplace-Beltrami
operator does not exceed $k_g\lambda$. As an immediate consequence
of Theorem \ref{thm-crit}, we get a version of this result for
functions from $\cF_{\lambda}$:

\begin{cor} \label{cor-Courant-2} For any $f\in\cF_\lambda$ and any
$a>0$, the number of nodal domains $A \in \cA $ with $m_A \geq a$
does not exceed $k_g\min\left(a^{-1}, a^{-2}\right)\lambda$.
\end{cor}

\noindent In particular, this applies to linear combinations of
the Laplace-Beltrami eigenfunctions on closed surfaces with the
eigenvalues $\leq \lambda$ (see Example \ref{ex-1} above). We
refer to \cite{A, GH} for a discussion on this subject.

Theorem \ref{thm-crit} is proved in Section \ref{sec-twoproofs}
below. Inequality \eqref{eq-gen-11} readily follows from the
Alexandrov-Backelman-Pucci inequality refined by Cabr\'{e}
\cite{C}, see inequality \eqref{eq-ABP} below. We take a different
route and obtain both inequalities in Theorem \ref{thm-crit} as a
special case of the following Kronrod-Yomdin-type estimate.

For a smooth function $f$ on $M$ and a regular value $c \in \R$ of
$f$, denote by $\beta(c,f)$ the number of connected components of
$f^{-1}(c)$. For a continuous function $u \in C(\R)$, define the
{\it generalized Banach indicatrix}
$$B(u,f) = \int_{-\infty}^{+\infty} u(c)\beta(c, f)\;dc\;.$$
In \cite{K}, Kronrod showed that the integral $\displaystyle
\int_{-\infty}^\infty \beta(c, f)\, dc$ is finite, and estimated
it from above through the {\em uniform norm} of the second
derivatives of $f$. This result was significantly extended in
different directions by many authors, notably, by Vitushkin,
Ivanov, and Yomdin, see \cite{Y} and references therein. In the
next theorem, we have in our disposal only $L^2$ bounds.

\begin{thm}\label{thm-main} For any function $ f\in\cF $ and
any continuous function $u$ on $\R$,
\begin{equation}\label{eq-ineq-0}
B(u,f) \leq k_g ||u\circ f|| \left( ||f|| + ||\Delta f|| \right)
\;.
\end{equation}
\end{thm}

\medskip
\noindent Theorem \ref{thm-main} is proved in Section
\ref{sec-Yomdin}.

Furthermore, we show that the nodal extrema of the eigenfunctions
of the Laplace-Beltrami operator satisfy the following additional
inequality. In the next theorem and the corollary we assume that
$M$ is a closed surface.

\begin{thm} \label{thm-eigen} For every Laplace-Beltrami
eigenfunction $f_{\lambda}$ with $||f_{\lambda}||=1$,
\begin{equation}\label{eq-gen-2}
\displaystyle \sum_{A\in \cA} m_A^6 \leq k_g \lambda^{3/2}\,.
\end{equation}
\end{thm}

\medskip\par\noindent Note that inequalities \eqref{eq-gen-11},\eqref{eq-gen-21}
and \eqref{eq-gen-2} for the Laplace-Beltrami eigenfunctions
$f_{\lambda}$ capture the sharp order of growth in $\lambda$. For
\eqref{eq-gen-11} and \eqref{eq-gen-21}, the example is the
sequence of eigenfunctions $f(x, y) = \sin nx \sin ny$,
$\lambda=n^2$, on the flat torus $\mathbb T^2$. For equation
\eqref{eq-gen-2}, the example is the sequence of zonal spherical
harmonics on the round sphere $\mathbb S^2$.

\medskip Curiously enough, inequality \eqref{eq-gen-2} complements
the classical bound \cite[Lemma 4.2.4]{Sogge} $$ \displaystyle
\max_{A \in \cA(f_{\lambda})} m_A = \displaystyle \max_M
|f_{\lambda}| \le k_g\lambda^{1/4}\;$$ as follows:

\begin{cor} \label{cor-Courant-1} Let
$f_{\lambda}$ be an eigenfunction of the Laplace-Beltrami operator
with the eigenvalue $\lambda$ so that $||f_{\lambda}||=1$. Then
for each $a>0$, the number of nodal domains $A$ of $f_{\lambda}$
with $m_A \geq a\lambda^{1/4}$ does not exceed $k_g a^{-6}$. In
particular, for fixed $a$, it remains bounded as
$\lambda\to\infty$.
\end{cor}

\medskip
\noindent Indeed, writing $n$ for the number of such nodal
domains, we get from \eqref{eq-gen-2} that
$$n(a\lambda^{1/4})^6 \leq \displaystyle \sum_{A\in \cA} m_A^6 \leq k_g
\lambda^{3/2}\,,$$ which yields the corollary.

It would be interesting to detect further restrictions on the
sequence of nodal maxima $\{m_A\}$ for the eigenfunctions. For the
analogue of Theorem \ref{thm-eigen} for the Dirichlet
eigenfunctions on surfaces with boundary, see Remark
\ref{rem.sogge} below.

Estimate \eqref{eq-gen-2} is a juxtaposition of the inner radius
theorem for nodal domains \cite{Mangoubi} with Sogge's $L^6$-bound
on the eigenfunctions, see Section~\ref{sec-Sogge} for the
details.

\section{Two approaches to nodal extrema}
\label{sec-twoproofs}

{\bf Proof of Theorem \ref{thm-crit}:} Take $f\in \cF_\lambda$ and
note that in this case the right hand side of \eqref{eq-ineq-0}
does not exceed $k_g\lambda$. For a nodal domain $A$ of $f$ and a
regular value $t$ of $f$ define $\beta_A(t)$ as the number of
connected components of $f^{-1}(t) \cap A$. Put $\epsilon_A =
\text{sign}(f\big|_A) \in \{-1;+1\}$. Observe that if the function
$u$ is even, then
\[
B(u,f) = \sum_{A \in \cA} \int_0^{m_A} u(t) \beta_A(\epsilon_A
t)dt\;.
\]
Since $\beta_A(\epsilon_A t) \geq 1$ almost everywhere on
$[0;m_A]$, we see that
\[
B(u,f) \geq \sum_{A \in \cA} \int_0^{m_A} u(t) \,dt\;.
\]
Choosing $u(t) = 1$ and applying Theorem \ref{thm-main}, we get
inequality \eqref{eq-gen-11}. To get \eqref{eq-gen-21}, we choose
$u(t) = |t|$. \Qed

\medskip

Another approach to inequality \eqref{eq-gen-11} is based on the
Alexandrov-Backelman-Pucci-Cabr\'{e} inequality \eqref{eq-ABP}
below. Let us illustrate the argument in the case when $M$ is a
plane domain equipped with a conformally Euclidean metric $g =
q(x,y)(dx^2+dy^2)$ with $0< q_- \leq q(x,y) < q_+$. Cabr\'{e}
\cite{C} showed that there exists a constant $C=C(q_-,q_+)$ such
that
\begin{equation}\label{eq-ABP}
\max_{M}|f| \leq C ( \text{Area}_g(M))^{1/2}||\Delta f||
\end{equation}
for every smooth function $f$ which vanishes on $\partial M$. It
is crucial that inequality \eqref{eq-ABP} remains valid with the
same constant $C$ when instead of $M$ we consider any domain $A
\subset M$. Given a function $f \in \cF_{\lambda}$ and a nodal
domain $A \in \cA(f)$ we get
\begin{equation} \label{eq-ABP-1}
m_A \leq C (\text{Area}_g(A))^{1/2} \Big{(} \int_A |\Delta
f|^2\;d\sigma \Big{)}^{1/2}\;, \end{equation} where $\sigma$
stands for the Riemannian measure on $M$. Let us sum up
inequalities \ref{eq-ABP-1} over all nodal domains $A \in \cA(f)$
and apply the Cauchy-Schwarz inequality. We get that
$$\sum_{A \in \cA(f)}m_A \leq C ( \text{Area}_g(M))^{1/2}||\Delta
f||\;,$$ which proves \eqref{eq-gen-11}.

\medskip

In order to extend this argument to general surfaces and
Riemannian metrics, one should use partition of unity associated
to a covering of $M$ by conformally Euclidean charts. We omit the
details.

Both approaches to inequality \eqref{eq-gen-11} presented above
have a similar geometric ingredient. The proof of inequality
\eqref{eq-ABP} involves, due to Alexandrov, the Gauss map of the
graph of the function $f$. Our proof of the Kronrod-Yomdin-type
bound, as we shall see in the next section, uses the Gauss map of
the level sets of $f$.

\section{Proof of the Kronrod-Yomdin-type \\ inequality} \label{sec-Yomdin}

Our proof of Theorem \ref{thm-main} is based on the following
strategy: Assume for simplicity that $u \equiv 1$. We shall find a
suitable ``length-type" functional $L$ so that its value on each
connected component of every regular level set of $f$ is greater
than $k_g
>0$. Thus its value $L(f^{-1}(c))$ on the full regular level
$f^{-1}(c)$ is at least $k_g\beta(c,f)$. Integrating against $dc$
we get that $\displaystyle \int_{-\infty}^\infty \beta(c, f) \, dc
\leq k_g^{-1}\int_{-\infty}^\infty L(c)\, dc$. Our choice of $L$
will enable us to rewrite the integral on the right hand side as
an integral over $M$ by using the co-area formula, and to estimate
it in terms of the $L^2$-norms of $f$ and its Laplacian. The usual
Riemannian length does not fit to the role of $L$ since $f$ may
have short level curves (e.g. in a neighborhood of a
non-degenerate maximum). However we observe that even short levels
become long when lifted to the unit circle bundle of $M$ together
with their normals (that is, via the Gauss map): indeed, the
normal field makes the full turn along such a curve. This will be
formalized below with the help of the Sasaki metric on the unit
circle bundle over $M$.

\medskip
The Levi-Civita connection gives us the canonical splitting $T(TM)
= \cV\oplus \cH$ into the vertical and the horizontal subspaces.
Each of them is canonically identified with $TM$. The metric $g
\oplus g$ on $TM$ (understood in the sense of the above splitting)
is called the {\it Sasaki metric}. Let
$$SM = \{(x,\xi) \in TM\;:\; |\xi|=1\}$$
be the unit circle bundle over $M$. Denote by $\rho$ the metric on
$SM$ induced from the Sasaki metric. Let $\kappa$ be the systole
of $(SM,\rho)$, that is $$\kappa = \inf
\text{length}_{\rho}(\alpha)\;,$$ where the infimum is taken over
all non-contractible closed curves $\alpha$ in $SM$.

\medskip Take any smooth function $f \in \cF$. Let $\gamma$
be a connected component of a regular level set of $f$. Let
$\nu=\nabla f/|\nabla f|$ be the field of unit normals along
$\gamma$. Choose the parameterization $\gamma(t)$ of $\gamma$ by
the Riemannian length so that is $|\dot{\gamma}| =1$. Denote by
$H_f\colon TM \to TM$ the Hessian of $f$. Put
$\widetilde{\gamma}(t) = (\gamma(t),\nu(t)) \in SM$. Thus
$\widetilde\gamma$ is the lift of $\gamma$ to $SM$. We start with
the following calculation:

\medskip
\noindent
\begin{lemma} \label{lem-geom} The length of the tangent vector to
$\widetilde\gamma$ with respect to the Sasaki metric is given by
$$|\dot{\widetilde\gamma}|^2_{\rho} =
1+\frac{(H_f\dot{\gamma},\dot{\gamma})^2}{|\nabla f|^2}\;.$$
\end{lemma}

\begin{proof} We write $\nabla$ for the covariant derivative with
respect to the Levi-Civita connection. Denote $w = \dot{\gamma}$.
Differentiating the identity $(\nu(t),\nu(t))=1$ we get that
$(\nabla_w \nu,\nu)=0$ and hence $\nabla_w\nu= (\nabla_w \nu,w)w$.
Using that
$$\nabla_w \nu =\nabla_w \frac{\nabla f}{|\nabla f|} =
\frac{1}{\nabla f} \nabla_w\nabla f + \Big{(} \nabla (|\nabla
f|^{-1}),w\Big{)}\nabla f$$ and that $H_fw = \nabla_w\nabla f$, we
get
$$\nabla_w \nu = \frac{(H_f w,w)}{|\nabla f|} w\;.$$
Using now that $|\dot{\widetilde\gamma}|^2_{\rho}= |w|^2 +
|\nabla_w \nu|^2$, we get the statement of the lemma.
\end{proof}

\medskip We shall need a simple (and well known) topological fact:
\begin{lemma} \label{lem-topo}
The curve $\widetilde\gamma$ is not contractible in $SM$.
\end{lemma}
\begin{proof} This is obvious when $\gamma$ is non-contractible in $M$.
Further, when $\gamma$ is contractible in $M$, the curve
$\widetilde\gamma$ is homotopic to the fiber, say $F_M$, of the
bundle $SM \to M$. The rest of the argument splits into three
cases.

\medskip
\noindent {\sc Case I:} $M=\mathbb{S}^2$. The manifold $SM$ equals
$\R P^3$. In this case the fiber $F_M$ represents the generator of
$\pi_1(\R P^3) = \Z_2$ and therefore is not contractible.

\medskip
\noindent {\sc Case II:} $M = \R P^2$. Look at the double cover
$p\colon S(\mathbb{S}^2) \to S(\R P^2)$. Observe that $p$ induces
a homeomorphism of the corresponding fibers $\widetilde{F}:=
F_{\mathbb{S}^2}$ and $F:= F_{\R P^2}$ . Assume on the contrary
that $F$ is contractible in $S(\R P^2)$. Applying the covering
homotopy theorem we get that $\widetilde{F}$ can be homotoped to
the fiber of $p$ in $S(\mathbb{S}^2)$. This contradicts to the
conclusion of Case I.

\medskip
\noindent {\sc Case III:} $\pi_2(M) = 0$. In this case the exact
homotopy sequence of the fibration $SM \to M$
$$0=\pi_2(M) \to \Z = \pi_1(F) \to \pi_1(SM)\;$$
shows that $F$ is not contractible in $SM$.
\end{proof}

\medskip
\noindent {\bf Proof of Theorem \ref{thm-main}:} For a regular
value $c$ of $f$, decompose the level set $f^{-1}(c)$ into
connected components
$$f^{-1}(c) = \gamma_1 \cup \dots \cup \gamma_{\beta(c,f)}\;.$$
Put $$L(c) =\sum_{i=1}^{\beta(c)} \text{length}_{\rho} {\widetilde
\gamma}_i\;.$$ Denote by $dt$ the length element along $f^{-1}(c)$
and by $d\sigma$ the Riemannian measure on $M$. By Lemma
\ref{lem-geom},
$$L(c) \leq \int_{f^{-1}(c)} (1+|H_f|^2/|\nabla f|^2)^{1/2}\; dt\;,$$
where $|H_f|$ stands for the operator norm of the Hessian. On the
other hand, applying Lemma \ref{lem-topo} we get that
$$L(c) \geq \kappa\beta(c,f)\;.$$
Thus
\begin{multline*}
B(u,f) \leq \kappa^{-1} \int_{-\infty}^{+\infty}
|u(c)|L(c)\;dc\; \\
\leq \kappa^{-1} \int_{-\infty}^{+\infty} dc\;
\int_{f^{-1}(c)}|u\circ f|\cdot(1+|H_f|^2/|\nabla f|^2)^{1/2}\;
dt\,.
\end{multline*}
Note that the iterated integral on the right-hand side equals
$\displaystyle \lim_{\epsilon\downarrow 0} I_\epsilon$ with
$$
I_\epsilon = \int_{-\infty}^{+\infty} dc\; \int_{f^{-1}(c)}|u\circ
f|\cdot \left( 1 + \frac{|H_f|^2}{\epsilon+|\nabla f|^2}
\right)^{1/2}\; dt\,.
$$
By the co-area formula
\begin{multline*}
I_\epsilon = \int_M |\nabla f|\cdot |u\circ f|\cdot \left( 1 +
\frac{|H_f|^2}{\epsilon+|\nabla f|^2} \right)^{1/2}\; d\sigma
\\ \le \int_M |u\circ f|\cdot \left( |\nabla f|^2 +
|H_f|^2\right)^{1/2}\; d\sigma\;.
\end{multline*}
Applying Cauchy-Schwarz inequality we see that this does not
exceed
\begin{equation}\label{eq-vilka-1} ||u\circ f|| \Big{(} \int_M (|\nabla
f|^2 + |H_f|^2)\; d\sigma \Big{)}^{1/2} \leq k_g||u\circ f||\cdot
\left( ||\nabla f|| + ||H_f|| \right)\;,
\end{equation}
where $||H_f||$ stands for the $L^2$-norm of $|H_f|$. The
classical a priori estimate \cite[Section~2.2]{Elliptic} (cf.
\cite[Chapter~II]{Lad}) tells us that
$$||H_f|| \leq k_g(||\Delta f||+ ||f||)$$
for every $f \in \cF$. Since $||\nabla f|| \le \frac12 (||f|| +
||\Delta f||)$, we are done. \Qed

\medskip
\noindent
\begin{rem}\label{rem-Bochner} {\rm
Assume that $f \in \cF$ is supported by the interior of $M$. For
instance, this holds automatically when $M$ has no boundary. In
this case we can modify the proof above and get the following
version of Theorem \ref{thm-main}:
\begin{equation}\label{eq-ineq-00}
B(u, f) \le k_{1} ||u\circ f|| \left( k_{2}||\Delta f||^2 +
k_{3}||\nabla f||^2\right)^{1/2} \,,
\end{equation}
where the constants $k_1,k_2,k_3$ have a transparent geometric
meaning. To have this estimate consistent in terms of units, we
introduce a real parameter $r>0$ which has the units of length and
consider a family of Sasaki metrics $\rho_r = r^2g\oplus g$ on
$SM$. We write $\kappa(r)$ for the systole of $(SM,\rho_r)$,
denote by $K$ the scalar curvature of $g$, and put $\displaystyle
K_{\rm min} = \min_M K(x)$. With this notation \eqref{eq-ineq-00}
holds with
\begin{eqnarray}\label{eq-ineq-1}
k_1 = \kappa^{-1}(r), \quad k_2 = r^2, \quad k_3 = 1-\frac12 r^2
K_{\rm min} \,.
\end{eqnarray}
For the proof of \eqref{eq-ineq-00} and \eqref{eq-ineq-1} we
repeat the arguments in the proof above and arrive at the estimate
\begin{equation}\label{eq-vsp-10}
B(u,f) \leq \kappa^{-1}(r)||u\circ f|| \Big{(} \int_M (|\nabla
f|^2 + r^2|H_f|^2)\; d\sigma \Big{)}^{1/2}\;,\end{equation} see
the left hand side of inequality \eqref{eq-vilka-1}. Instead of
proceeding as in \eqref{eq-vilka-1}, we use the
Bochner-Lichnerowicz formula (see \cite[p.85]{Ch}):
$$
\frac{1}{2}\Delta(|\nabla f|^2) = \operatorname{tr}(H_f^2) +
(\nabla f,\nabla \Delta f)+ \frac{1}{2}K |\nabla f|^2\;.
$$
Integrating it over $M$ we get that
\begin{equation}\label{eq-2}
\int_M \operatorname{tr}(H_f^2)\; d\sigma = ||\Delta f||^2 -
\frac{1}{2}\int_M K|\nabla f|^2 \;d\sigma \;.
\end{equation}
Combining this with \eqref{eq-vsp-10} and using that $|H_f|^2\le
\operatorname{tr}(H_f^2)$, we get the desired result. }
\end{rem}

\medskip
\noindent
\begin{example} \label{ex-disc}{\rm As an illustration, consider
the case when $(M,g)$ is any bounded simply connected domain with
smooth boundary in the Euclidean plane $\R^2(x,y)$. In this case
$SM = \mathbb{S}^1 \times M$, where the circle $\mathbb{S}^1$ is
identified with $\R/2\pi\Z$ and is equipped with the coordinate
$\theta$ $(\text{mod}\;2\pi)$. The Sasaki metric $\rho_r$ is given
by $r^2d\theta^2+dx^2+dy^2$. The systole $\kappa(r)$ is equal to
the length of the fiber $\mathbb{S}^1 \times \text{point}$, and so
$k_1 = \frac{1}{2\pi r}$. The curvature $K$ vanishes and hence
$k_3= 1$. Applying inequality \eqref{eq-ineq-00} and passing to
the limit as $r \to \infty$ we get that
$$B(u,f) \leq \frac{1}{2\pi}||u \circ f|| \cdot ||\Delta f||\;.$$
Taking $u=1$, we get
$$\max_M f - \min_M f \leq B(1,f) \leq B(1,f) \leq
\frac{1}{2\pi} (\text{Area}(M))^{1/2}||\Delta f||\;,$$ which is a
special case of the Alexandrov-Bakelman-Pucci-Cabr\'{e} inequality
\eqref{eq-ABP}.}
\end{example}

\medskip
\begin{rem}\label{rem.Gr}
{\rm As a by-product of our method, we get the following
inequality: Let $M\subset \R^2$ be a bounded plane domain with
smooth boundary equipped with the Euclidean metric. Let $f$ be a
smooth function vanishing on $\partial M$. Then
\begin{equation}\label{eq.Gr}
\max_M |f| \leq \frac{1}{2\pi} \int_M |H_f| d\sigma\;.
\end{equation}
This inequality is sharp: take $M$ to be the unit disc and $f(x,
y)=1-(x^2+y^2)$. Note that in contrast to inequality
\eqref{eq-ABP}, it involves the $L^1$-norm of the second
derivatives of $f$.

\medskip
To prove inequality \eqref{eq.Gr}, it suffices to show that
\begin{equation}\label{eq.1}
B(1,f) \leq \frac{1}{2\pi} \int_M |H_f| d\sigma\;.
\end{equation}
Introduce the Sasaki metric $\rho_r$ on $SM$ with parameter $r$ as
in Remark~\ref{rem-Bochner}. Note that the Sasaki length of the
lift to $SM$ of the field of normals along any simple closed curve
in $M$ is $\geq 2\pi r$. Arguing as in the proof of
Theorem~\ref{thm-main}, we get
$$
B(1,f) \leq (2\pi r)^{-1} \int_{-\infty}^{+\infty} dc\;
\int_{f^{-1}(c)}(1+r^2|H_f|^2/|\nabla f|^2)^{1/2}\; dt\,.$$
Furthermore,
$$(1+r^2|H_f|^2/|\nabla f|^2)^{1/2} \leq 1+ r|H_f|/|\nabla f|\;.
$$
Thus by the co-area formula
$$B(1,f) \leq (2\pi r)^{-1} \int_M (|\nabla f| +
r|H_f|)d\sigma\;.$$ This holds for every $r$. Passing to the limit
as $r \to +\infty$ we get the desired inequality \eqref{eq.1}.
\hfill $\Box$

\medskip
In fact, formally speaking, the Sasaki metric is not needed in the
proof of \eqref{eq.Gr}. Instead of dealing with the
$\rho_r$-length of $f^{-1}(c)$ and passing to the limit as $r \to
+\infty$, one can start with the quantity
$$\int_{-\infty}^{+\infty} dc\; \int_{f^{-1}(c)} |Q(t)|dt \;,$$
where $Q(t)$ is the curvature of $f^{-1}(c)$. }
\end{rem}

\begin{rem}{\rm It is not difficult to see that the equivalent way
to state Theorem~\ref{thm-main} is as follows:
\[
\int_{-\infty}^\infty \frac{\beta^2(c, f)}{l(c, f)} \, dc \le k_g
\left( ||f|| + ||\Delta f|| \right) \,,
\]
where
\[
l(c, f) = \int_{\{f=c\}} \frac{ds}{|\nabla f|} = \lim_{\epsilon\to
0 } \frac{\sigma\left\{|f-c|<\epsilon\right\}}{\epsilon}
\]
is the Leray length of the level line $\{f=c\}$. Here $\sigma$
stands for the Riemannian measure on $M$. We skip the details. }
\end{rem}

\section{More on nodal extrema of eigenfunctions}
\label{sec-Sogge}

In this section, we prove Theorem \ref{thm-eigen}. Consider a
nodal domain $A_i$, and take the point $p_i\in A_i$ such that
$|f_\lambda (p_i)| = m_i$. Without loss of generality, assume that
$f_\lambda$ is positive in $A_i$.

We choose local coordinates $(x,y)$, $ x^2 +y^2 \leq 2 $, near
$p_i$ so that the metric $g$ in these coordinates is conformally
Euclidean: $g = q(x,y)(dx^2 + dy^2)$. The point $p_i$ corresponds
to the origin: $p_i = (0,0)$. This choice can be made in such a
way that $ q(x,y) \leq k_g$. By disk $D(p_i, r)$ centered at a
point $p_i$ with radius $r$ we mean the set $\{x^2+y^2\le r^2\}$,
where $(x,y)$ are local conformal coordinates near $p_i$.

By the inradius theorem \cite[Lemma~10]{Mangoubi}, there exists a
value $\mu$ depending only on the metric $g$ such that the disk
$\displaystyle D(p_i, \mu \lambda^{-\frac12})$ is contained in
$A_i$. We also assume that $\mu\lambda^{-\frac12} \le 1$,
otherwise, we just choose a smaller $\mu$. Set
$r=\mu\lambda^{-\frac12}$. In local conformal coordinates, the
eigenfunction $f_{\lambda}$ satisfies the equation $\Delta_e
f_\lambda + \lambda q(x,y) f_\lambda = 0$, where $\Delta_e$ is the
Euclidean Laplacian. Representing the function $f_\lambda$ as a
sum of Green's potential and the Poisson integral in
$\{\sqrt{x^2+y^2}\le tr\}$ with $\frac12 \le t \le 1$, we get
\[
m_i = f_\lambda(0) = \frac{\lambda}{2\pi} \int_{|\zeta|<tr}
q(\zeta) f_\lambda(\zeta) \log\frac{tr}{|\zeta|}\,
d\sigma_e(\zeta) + \int_{|\zeta|=tr} f_\lambda(\zeta)\,
d\ell(\zeta)\,.
\]
Here, $\sigma_e$ is the Euclidean measure $dxdy$ and $\ell$ is the
Lebesgue measure on the circle $\{|\zeta|=tr\}$ which is
normalized so that the total measure of the circle equals $1$.

Using H\"{o}lder's inequality, we get
\[
m_i^6 \le k_g \left( \lambda \int_{|\zeta|<tr} f_\lambda^6\,
d\sigma_e + \int_{|\zeta|=tr} f_\lambda^6 \, d\ell \right)\,.
\]
Integrating this estimate with respect to $t$ from $\frac12 $ to
$1$, and taking into account that the Euclidean measure $\sigma_e$
and the Riemannian measure $\sigma$ are equivalent on the disks
$D(p_i, r)$, we get
\[
m_i^6 \le k_g \lambda \int_{D(p_i, r)} f_\lambda^6\, d\sigma\,.
\]
Note that the disks $D(p_i, r)$ do not overlap. Hence, summing up
by $i$, we get
\[
\sum m_i^6 \le k_g \lambda \int_X f_\lambda^6\, d\sigma\,.
\]
Classical Sogge's bound $||f_\lambda||_6 \le k_g
\lambda^{\frac{1}{12}}$ \cite[Chapter~5]{Sogge} completes the
proof of \eqref{eq-gen-2}.

\hfill $\Box$

\medskip
\begin{rem}\label{rem.sogge}
{\rm The inradius theorem can be extended to nodal domains of the
Dirichlet eigenfunctions on surfaces with boundary. Applying the
results of a recent paper by Smith and Sogge \cite{SS}, one can
show by the same argument as above that $$ \sum_{A\in\mathcal
A(f_\lambda)} m_A^8 \le k_g \lambda^2 $$ for any Dirichlet
eigenfunction $f_\lambda$. }
\end{rem}

\medskip
\noindent {\bf Acknowledgement.} We thank Fima Gluskin, Misha
Katz, Fedya Nazarov, Iosif Poltero\-vich, Zeev Rudnick and Yosi
Yomdin for many valuable discussions. We thank Daniel Grieser for
stimulating correspondence that yielded the second proof of
inequality \eqref{eq-gen-11} and Remark~\ref{rem.Gr}. We are
grateful to Kari Astala and Yehuda Pinchover for referring us to
the Alexandrov-Bakelman-Pucci inequality.

\end{document}